\documentclass[a4paper,11pt,oneside]{article}
\usepackage{amsmath,amsthm,amssymb,mathdots,hyperref,comment}
\usepackage{color}
\usepackage[normalem]{ulem}

\newcommand{\Z}{\mathbb{Z}}
\newcommand{\pres}[2]{\langle {#1}\ |\ {#2} \rangle}
\newtheorem{theorem}{Theorem}
\newtheorem{lemma}[theorem]{Lemma}
\newtheorem{corollary}[theorem]{Corollary}
\newtheorem{remark}[theorem]{Remark}
\newtheorem{maintheorem}{Theorem}

\newtheorem{maincorollary}[maintheorem]{Corollary}

\numberwithin{theorem}{section}
\numberwithin{equation}{section}

\allowdisplaybreaks

\DeclareMathOperator{\rank}{rank}

\DeclareMathOperator{\res}{Res}

\begin{document}
\title{Matrices in companion rings, Smith forms, and the homology of $3$-dimensional Brieskorn manifolds}
\author{Vanni Noferini\thanks{Corresponding author. \textsc{Department of Mathematics and Systems Analysis, Aalto University, P.O. Box 11000 (Otakaari 24), FI-00076 AALTO, Finland.} \textit{E-mail address}, \texttt{vanni.noferini@aalto.fi}} and Gerald Williams\thanks{ \textsc{Department of Mathematical Sciences, University of Essex, Wivenhoe Park, Colchester, Essex CO4 3SQ, UK.}
  \textit{E-mail address}, \texttt{Gerald.Williams@essex.ac.uk} }}

\maketitle

\begin{abstract}
We study the Smith forms of matrices of the form $f(C_g)$ where $f(t),g(t)\in R[t]$, where $R$ is an elementary divisor domain and $C_g$ is the companion matrix of the (monic) polynomial $g(t)$. Prominent examples of such matrices are circulant matrices, skew-circulant matrices, and triangular Toeplitz matrices. In particular, we reduce the calculation of the Smith form of the matrix $f(C_g)$ to that of the matrix $F(C_G)$, where $F,G$ are quotients of $f(t),g(t)$ by some common divisor. This allows us to express the last non-zero determinantal divisor of $f(C_g)$ as a resultant. A key tool is the observation that a matrix ring generated by $C_g$ -- the companion ring of $g(t)$ -- is isomorphic to the polynomial ring $Q_g=R[t]/<g(t)>$. We relate several features of the Smith form of $f(C_g)$ to the properties of the polynomial $g(t)$ and the equivalence classes $[f(t)]\in Q_g$. As an application we let $f(t)$ be the Alexander polynomial of a torus knot and $g(t)=t^n-1$, and calculate the Smith form of the circulant matrix $f(C_g)$. By appealing to results concerning cyclic branched covers of knots and cyclically presented groups, this provides the homology of all Brieskorn manifolds $M(r,s,n)$ where $r,s$ are coprime.
\end{abstract}


\section{Introduction}\label{sec:intro}

Let $R$ be a commutative ring (with unity) other than the trivial ring, fix a monic polynomial $g(t) = t^n + \sum_{k=0}^{n-1} g_k t^k \in R[t]$, and let $C_g$ be the companion matrix of $g(t)$. For $n\geq 2$, the subset of $R^{n\times n}$ consisting of matrices $f(C_g)$ that are polynomials in $C_g$ with coefficients in $R$ forms a commutative ring, which we call the \em companion ring \em of $g(t)$ and denote by $R_g$. Important and well studied rings of matrices arise as special cases: if $g(t)=t^n$, then $R_g$ is the commutative ring of lower triangular $n \times n$ Toeplitz matrices \cite{BP94} with entries in $R$; if $g(t)=t^n-1$, then $R_g$ is the commutative ring of $n \times n$ circulant matrices \cite{Davis} with entries in $R$; if $g(t)=t^n+1$, then $R_g$ is the commutative ring of $n \times n$ skew-circulant matrices \cite{Davis} with entries in $R$.

When $R$ is an integral domain, $g(t)$ has $n$ roots (counted with multiplicities) in some appropriate extension of $R$ and, for $f(t)\in R[t]$, the determinant of $f(C_g)$ can be expressed as the resultant
\begin{equation}\label{eq:detfCg}\det f(C_g) = \prod_{\theta : g(\theta) = 0} f(\theta) =: \res(f,g).\end{equation}
Note that in the last equality we are implicitly fixing the normalization in the definition of the resultant; we will follow this choice throughout. (See Section~\ref{sec:prelims} for definitions of undefined terms and notation used in this Introduction, together with relevant background). Restricting to the case that $R$ is an elementary divisor domain, such as the ring of integers, one may seek to study the Smith forms of matrices $f(C_g)$. Our first main result shows that $f(C_g)$ is equivalent to the direct sum of a matrix $F(C_G)$ and a zero matrix, where $F,G$ are quotients of $f,g$ by any of their (monic) common divisors, and so relates the Smith forms of $f(C_g)$ and $F(C_G)$.

\begin{maintheorem}\label{thm:C}
Let $g(t) \in R[t]$ be monic of degree $n$, and let $f(t) \in R[t]$ where $R$ is an elementary divisor domain. Suppose that $g(t) = G(t) z(t)$, $f(t) = F(t) z(t)$ where $z(t)$ is a monic common divisor of $f(t)$ and $g(t)$. Then  $f(C_g) \sim F(C_G) \oplus 0_{m\times m}$, where $m= \deg z(t)$. In particular, $F(C_G)$ has invariant factors $s_1,\dots,s_r$ if and only if $f(C_g)$ has invariant factors $s_1, \dots, s_r, 0$ (repeated $m$ times).
\end{maintheorem}

An immediate corollary (Corollary~\ref{cor:D}) expresses the last non-zero determinantal divisor as the resultant of $F(t)$ and $G(t)$. This therefore generalizes the expression~(\ref{eq:detfCg}) to the case of singular matrices $f(C_g)$.

\begin{maincorollary}\label{cor:D}
In the notation of Theorem~\ref{thm:C}, suppose that $z(t)$ is the monic \emph{greatest} common divisor of $f(t)=z(t)F(t)$ and $g(t)=z(t) G(t)$. Then the last non-zero determinantal divisor of $f(C_g)$ is (up to units of $R$)
\begin{equation*}\label{eq:gamma_r}\gamma_r = \prod_{\theta : G(\theta) = 0} F(\theta) = \res(F,G).\end{equation*}
\end{maincorollary}

As an application, in Theorem~\ref{thm:dunwoody}, we calculate the Smith form of the integer matrix $f(C_g)$ where $f(t)$ is the Alexander polynomial of the torus knot $K(r,s)$, i.e.
\begin{alignat}{1}
f(t)=\frac{(t^{rs}-1)(t-1)}{(t^s-1)(t^r-1)}\label{eq:AlexanderTorus}
\end{alignat}
 and $g(t)=t^n-1$. As we explain in Section~\ref{subsec:cycpresgps} this allows us, in Corollary~\ref{cor:HomologyDunwoody}, to calculate the homology of all $3$-dimensional Brieskorn manifolds $M(r,s,n)$ where $r,s$ are coprime. This generalizes (part of)~\cite[Proposition~5]{CavicchioliAnnali}, which deals with the case $r=2$.

\begin{maintheorem}\label{thm:dunwoody}
Let $r,s$ be coprime positive integers, $n\geq 2$,  such that $x:=(r,n)\leq y:=(s,n)$ and let $f(t)\in \Z[t]$ be the Alexander polynomial of the torus knot $K(r,s)$ as in~(\ref{eq:AlexanderTorus}), $g(t)=t^n-1\in \Z[t]$.
Then the Smith form of $f(C_g)$ has non-unit invariant factors: $\frac{r}{x}$ (repeated $y-x$ times); $\frac{rs}{xy}$ (repeated $x-1$ times); $0$ (repeated $(x-1)(y-1)$ times).
\end{maintheorem}

We note that there is no loss of generality in assuming, as in the statement of Theorem \ref{thm:dunwoody}, that $(r,n) \leq (s,n)$ for if not we may simply swap the roles of $r$ and~$s$.

\begin{maincorollary}\label{cor:HomologyDunwoody}
Let $r,s,n\geq 2$ where $r$ and $s$ are coprime. Then setting $x:=(r,n)$ and $y:=(s,n)$, the homology of the $3$-dimensional Brieskorn manifold  $M=M(r,s,n)$ is
\[ H_1(M)\cong \begin{cases}
  \Z_{r/x}^{y-x} \oplus \Z_{rs/(xy)}^{x-1}\oplus \Z^{(x-1)(y-1)} &\mathrm{if}~x\leq y,\\
  \Z_{s/y}^{x-y} \oplus \Z_{rs/(xy)}^{y-1}\oplus \Z^{(x-1)(y-1)} &\mathrm{if}~y\leq x.
\end{cases}
\]
\end{maincorollary}

\section{Preliminaries}\label{sec:prelims}

\subsection{Smith Forms and Elementary Divisor Domains}\label{subsec:SmithFormsEDD}

Given a GCD domain $R$, we denote the greatest common divisor\footnote{Although, strictly speaking, the greatest common divisor is only defined up to multiplication by units of $R$, we assume here and throughout that a choice is made by some arbitrary, but fixed, choice of normalization. To give a concrete example, for $R=\Z$ one may choose gcds to be always non-negative integers.} of the $n$-tuple $a_1,\dots,a_n \in R$ by $(a_1,\dots,a_n)$. An \emph{elementary divisor domain} (EDD) \cite[p.\,16]{Friedland} is  an integral domain $R$ such that, for any triple of elements $a,b,c \in R$, there exist $x,y,z,w \in R$ satisfying $(a,b,c) = zxa+zyb+wyc.$

By choosing $c=0$ in this definition, it follows that every EDD is a B\'{e}zout domain; it is conjectured, but to our knowledge still an open problem, that the converse is false \cite{K49,L12}. Every principal ideal domain (PID) is an EDD see, for example, \cite[Theorem 1.5.3]{Friedland}; a classical example of an EDD that is not a PID is the ring of functions that are holomorphic on a simply connected domain \cite{H43,Q97}.

The following classical theorem is named after H. J. S. Smith, who studied the case $R=\mathbb{Z}$~\cite{Smith}. Frobenius proved the Smith Theorem in~\cite{Frobenius} assuming that $R$ is a ring of univariate polynomials with coefficients in a field. For a proof of the theorem when $R$ is an EDD (the weakest permissible assumption under which it can hold) see, for example, \cite[Theorem 1.14.1]{Friedland}.

To state the Smith Theorem, recall \cite[p.\,12]{Newmanbook} that a square matrix $U
\in R^{n \times n}$ is called \emph{unimodular} if $
\det U$ is a unit of the base commutative ring $R$. Equivalently, unimodular matrices are precisely those matrices that are invertible over $R$, i.e., whose inverse exists and also belongs to $R^{n \times n}$.

\begin{theorem}[Smith Theorem]\label{thm:SNF}
Let $R$ be an EDD and $M \in R^{m \times n}$. Then there exist unimodular matrices $U \in R^{n \times n}$, $V \in R^{m \times m}$ such that $U M V = S$ where $S$ is diagonal and satisfies $S_{i,i} \mid S_{i+1,i+1}$ for all $i=1,\dots,\min(m,n)-1$. Further, let $\gamma_0=1 \in R$, and for $i=1,\dots,\min(m,n)$ define the $i$-th  \emph{determinantal divisor} $\gamma_i$ to be the greatest common divisor (gcd) of all minors of $M$ of order $i$. Then
$$S_{i,i} = \frac{\gamma_i}{\gamma_{i-1}} =: s_i(M),$$
where the diagonal elements $s_i(M)$, $i=1,\dots,\min(m,n)$, are called the \emph{invariant factors} of $M$. The matrix $S$ is called the
\emph{Smith form} of $M$, and it  is uniquely determined by $M$ up to multiplication of the invariant factors by units of $R$.
\end{theorem}

To make the Smith form $S$ uniquely determined by $M$, one might consider an appropriate normalization of the determinantal divisors, or equivalently of the invariant factors. This  ``appropriate normalization''  is a conventional, albeit arbitrary, choice that depends on the base ring $R$. For instance, typical requirements are that the invariant factors are non-negative integers when $R=\mathbb{Z}$; or that the invariant factors are monic polynomials when $R=\mathbb{F}[x]$ (univariate polynomials with coefficients in a field $\mathbb{F}$).  To avoid pedantic repetitions of sentences like ``up to multiplication by units of $R$'', we assume that one such normalization is tacitly agreed on.

Generally, if there are two unimodular $U \in R^{m \times m}, V \in R^{n \times n}$ such that $UMV=N$, we say that $M$ and $N$ are \emph{equivalent} (over $R$) and write $M \sim N$. If furthermore $V=U^{-1}$ then $M$ and $N$ are said to be \emph{similar} over $R$. The Smith Theorem can therefore be stated as follows: every matrix with entries in an EDD is equivalent to a diagonal matrix whose diagonal elements form a divisor chain. This  is, in fact, a characterization of an EDD, usually taken as a definition \cite{K49,L12,Q97}: an EDD is an integral domain $R$ over which the Smith Theorem holds. We mention two immediate consequences of the Smith Theorem that will be useful to us: firstly, any pair of $m \times n$ matrices with entries in $R$ are equivalent if and only if they have the same invariant factors; secondly, since rank is preserved by multiplication by invertible matrices, $M$ has rank $r$ if and only if its invariant factors satisfy $s_i(M)=0$ precisely when $i>r$.

\subsection{Cyclically presented groups and Brieskorn manifolds}\label{subsec:cycpresgps}

Any finitely generated abelian group $A$ is isomorphic to a group of the form $A_0\oplus \Z^\beta$ where $A_0$ is a finite abelian group and $\beta\geq 0$. The number $\beta=\beta(A)$ is called the \em Betti number \em (or \em torsion-free rank\em) of $A$. Clearly $A$ is infinite if and only if $\beta(A)\geq 1$ and $A$ is a free abelian group if and only if $A_0=1$.

Given a group presentation $P=\pres{x_0,\ldots,x_{n-1}}{R_0,\ldots ,R_{m-1}}$ ($n,m\geq 1$) the \em relation matrix \em of $P$ is the
$n\times m$ integer matrix $M$ whose $(i,j)$ entry is the exponent sum of generator $x_i$ in relator $R_j$. If the rank of $M$ is $r$ and the invariant factors of the Smith Form of $M$ are $s_1,\ldots,s_{n-r}$ then the abelianization of the group $G$ defined by the presentation $P$ is
\[ G^\mathrm{ab}\cong \Z_{s_1}\oplus \ldots \oplus \Z_{s_{r}}\oplus \Z^{n-r}.\]
(See, for example, \cite[pp.\,146--149, Theorem~3.6]{MKS} or \cite[pp.\,54--57, Theorem~5]{JohnsonBook}.) Thus $\beta(G^\mathrm{ab})=n-r$ and if $G^\mathrm{ab}=A_0\oplus \Z^\beta$ then $|A_0|=|\prod_{i=1}^r s_i|$, i.e. the last non-zero determinantal divisor $\gamma_r$ of $M$.

A \em cyclic presentation \em is a group presentation of the form
\begin{alignat*}{1}
P_n(w)=\pres{x_0,\ldots ,x_{n-1}}{w(x_i,x_{i+1},\ldots ,x_{i+n-1})\ (0\leq i<n)}
\end{alignat*}
where $w=w(x_0,x_1,\ldots ,x_{n-1})$ is some fixed element of the free group $F(x_0,\ldots ,x_{n-1})$ and the subscripts are taken mod~$n$, and the group $G_n(w)$ it defines is called a \em cyclically presented group\em. If, for each $0\leq i<n$, the exponent sum of $x_i$ in $w(x_0,\ldots ,x_{n-1})$ is $a_i$ then the relation matrix of $P_n(w)$ is the circulant matrix $C$ whose first row is $(a_0,a_1,\ldots, a_{n-1})$. The \em representer polynomial \em of $C$ is the polynomial
\[ f(t)=f_w(t)=\sum_{i=0}^{n-1}a_it^i \in \Z[t]\]
and setting $g(t)=t^n-1\in \Z[t]$, $R=\Z$, the relation matrix of $P_n(w)$ is the circulant matrix $f_w(C_g)$. Thus, results concerning the Smith forms of such matrices $f(C_g)$ provide information about the abelianization of the cyclically presented group $G_n(w)$.

This, in turn, allows us to calculate the homology of certain $3$-dimensional manifolds, as we now describe. For a $3$-manifold $M$, the first homology $H_1(M)$ is isomorphic to the abelianization of its fundamental group (see, for example, \cite[Theorem~2A.1]{Hatcher}). Thus, given a $3$-manifold whose fundamental group has a cyclic presentation $P_n(w)$ with representer polynomial $f_w(t)$, the Smith form of the integer circulant $f_w(C_g)$ provides the homology of~$M$. Suitable manifolds include, for example, all Dunwoody manifolds~\cite{Dunwoody}.

\em Brieskorn manifolds \em were introduced in~\cite{Brieskorn} and the $3$-dimensional Brieskorn manifolds $M(r,s,n)$ ($r,s,n\geq 2$) were studied by Milnor in~\cite{Milnor}. As noted in~\cite[p.\,176]{Milnor}, the order of the homology of $M(r,s,n)$ was computed by Brieskorn~\cite{Brieskorn}, who showed that the homology is trivial if and only if $r,s,n$ are pairwise relatively prime. An algorithm for computing the homology itself, conjectured in \cite{Orlik}, was proved in~\cite{Randell}. Further, the homology was calculated for the case $r=2$ in~\cite{CavicchioliAnnali}. The manifolds $M(r,s,n)$ can be described as $n$-fold cyclic branched covers of the $3$-sphere $S^3$ branched over the torus link $K(r,s)$, or torus knot $K(r,s)$ when $(r,s)=1$~\cite[Lemma~1.1]{Milnor}. Torus knots lie in a very general class of knots called \em $(1,1)$-knots. \em A special case of~\cite[Theorem~3.1]{MulazzaniJKorea} is that if a manifold $M$ is an $n$-fold cyclic branched cover of $S^3$ branched over a $(1,1)$-knot then its fundamental group $\pi_1(M)$ has a cyclic presentation $P_n(w)$ for some $w$. Moreover, by~\cite[Proposition~7]{CattagbrigaMulazzaniMPCPS} (see also~\cite[Theorem~4]{CattabrigaAlexanderPolynomial})  $w$ can be chosen so that the representer polynomial $f_w(t)$ of $P_n(w)$ is equal to the projection of the Alexander polynomial $\Delta_K(t)$ of $K$ to $\Z[t]/<t^n-1>$. Thus, in particular, if $(r,s)=1$ then the $3$-dimensional Brieskorn manifold $M(r,s,n)$ has a cyclic presentation $P_n(w)$ where $f_w(t)$ is equal to the projection of the Alexander polynomial $f(t)$ given at~(\ref{eq:AlexanderTorus}) to $\Z[t]/<g(t)>$ where $g(t)=t^n-1$. Hence, for coprime $r,s$, the calculation of the Smith form of $f(C_g)$ given in Theorem~\ref{thm:dunwoody} provides the homology of $M(r,s,n)$, as in Corollary~\ref{cor:HomologyDunwoody}.

\section{Quotient polynomial rings as a ring of matrices}\label{sec:quotientpolyringasmatrices}

As in the Introduction, let $R$ be a commutative ring (with unity) other than the trivial ring and fix a monic polynomial $g(t) = t^n + \sum_{k=0}^{n-1} g_k t^k \in R[t].$ The \emph{quotient ring} $Q_g = R[t]/\langle g(t) \rangle$ is the ring of the equivalence classes of polynomials in $R[t]$ modulo the ideal generated by $g(t)$. Specifically, given any $f(t) \in R[t]$, one defines the equivalence class $[f(t)] \in Q_g$ as
$$[f(t)]:= \{ h(t) \in R[t] : h(t) \equiv f(t)~\mathrm{mod}~g(t)   \};$$ here and below, $h(t) \equiv f(t)$~mod~$g(t)$ is a compact notation to mean that there exists $q(t) \in R[t]$ such that $f(t) = h(t) + q(t) g(t)$. On the other hand, associated with $g(t)$ is its \emph{companion matrix}
$$ C_g = \begin{bmatrix} -g_{n-1} & \dots & -g_1 & - g_0\\
1 & & & \\
& \ddots & &\\
& & 1 &
\end{bmatrix} \in R^{n \times n},$$
(where, as throughout this paper, entries not explicitly displayed are assumed to be $0$); observe that
the matrix $C_g$ is a representation of the multiplication-by-$[t]$ operator in the quotient ring $Q_g$, in that
\begin{equation}\label{eq:CgVander=tVander}
C_g \begin{bmatrix} t^{n-1} & \cdots & t& 1 \end{bmatrix}^T \equiv t \begin{bmatrix} t^{n-1} & \cdots & t& 1 \end{bmatrix}^T ~\mathrm{mod}~g(t).
\end{equation}
For more details on this viewpoint on companion matrices, as well as some generalizations, see, for example, \cite[Section~2]{NN16}, \cite[Section~9]{NP15},  and the references therein. It is well known in matrix theory that the characteristic polynomial of $C_g$ is precisely $g(t)$: see, for example, the proof of \cite[Theorem 1.1]{GohLR09}; although there it is assumed that $R=\mathbb{C}$, the proof is purely algebraic and only requires that $R$ is a commutative ring.

Theorem~\ref{thm:zero} below is a special case of the First Isomorphism Theorem for rings, but we will give an elementary matrix theoretical proof. It illustrates how to expand the idea of mapping $[t]$ to $C_g$, and why it generates a matrix algebra. Given an equivalence class $[f(t)] \in Q_g$, it is natural to consider the coefficients $f_k \in R$ of its expansion in the monomial basis of $Q_g$, i.e.,
$$[f(t)] = \sum_{k=0}^{n-1} f_k [t^k].$$
Note that the above notation is generally equivalent to $f(t) \equiv$ \linebreak $\sum_{k=0}^{n-1} f_k t^k$~mod~$g(t)$; the actual equality $f(t) = \sum_{k=0}^{n-1} f_k t^k$ holds if and only if $\deg f(t)<n$. Before stating Theorem~\ref{thm:zero} we recall that the Cayley-Hamilton theorem holds for matrices over any commutative ring \cite[Problem 2.4.P3]{HJ85}. Thus, if $f(t) \equiv h(t)$~mod~$g(t)$ then for some $q(t) \in R[t]$ we have $f(C_g)=h(C_g) + g(C_g) q(C_g) = h(C_g)$. Hence, writing $\sum_{k=0}^{n-1} f_k C_g^k=f(C_g)$ is consistent even if, in general, one may be taking $\deg f(t) \geq n$ so that the coefficients of $f(t)$ differ from the coefficients of $[f(t)]$.

It is easy to show that, under the usual matrix addition and multiplication, the subset of $R^{n \times n}$ consisting of matrices that are polynomials in $C_g$ with coefficients in $R$ forms a commutative ring,  which we will denote by $R_g$. Theorem~\ref{thm:zero} shows that this is isomorphic to the quotient ring $Q_g$, and is fundamental to our methods.

\begin{theorem}\label{thm:zero}
The map $\mathcal{M}: Q_g \rightarrow R_g$ given by
\[[f(t)]= \sum_{k=0}^{n-1} f_k [t^k] \mapsto \mathcal{M}([f(t)]) = \sum_{k=0}^{n-1} f_k C_g^k = f(C_g)\]
is a ring isomorphism.
\end{theorem}

\begin{proof}
The map $\mathcal{M}$ is clearly bijective, maps $[0]$ to $0$ and $[1]$ to $I$, and satisfies $\mathcal{M}([f_1(t)]+[f_2(t)]) = \mathcal{M}([f_1(t)]) + \mathcal{M}([f_2(t)])$. If $f_1(t) f_2(t) = f_3 (t) + q(t) g(t)$, $\deg f_3(t) < n$, then using the Cayley-Hamilton theorem
\begin{alignat*}{1}
\mathcal{M}([f_1(t)]) \mathcal{M}([f_2(t)])&= f_1(C_g) f_2(C_g)= f_3(C_g) + q(C_g) g(C_g)\\
&\qquad = f_3(C_g)= \mathcal{M}([f_3(t)])= \mathcal{M}([f_1(t) f_2(t)]).
\end{alignat*}
\end{proof}

Theorem~\ref{thm:zero} shows that, for every monic polynomial $g(t) \in R[t]$,  we can define a matrix algebra $R_g$ that satisfies the important property of being isomorphic to the quotient ring $Q_g$. It is also useful to observe that, by \eqref{eq:CgVander=tVander}, for all $j=0,\dots,n-1$ the $(n-j)$-th row of $f(C_g)$ contains the coefficients of $[t^j f(t)]$ in the monomial basis $[1],\dots,[t^{n-1}]$. In particular, the last row of $f(C_g)$ is precisely
$$\begin{bmatrix} f_{n-1} & \dots & f_1 & f_0\end{bmatrix}.$$
As mentioned in the Introduction, the commutative rings of lower triangular $n \times n$ Toeplitz matrices, of circulant matrices, and of skew-circulant matrices all arise as special cases of $R_g$.

If we now specialize to the case where $R$ is an EDD then, given $g(t) \in R[t]$ (monic) and $[f(t)] \in Q_g$, it makes sense to study the Smith canonical form of $f(C_g)$. An important example of an EDD is the ring of the integers~$\mathbb{Z}$ and in this setting, $g(t)$ is a monic integer polynomial, $[f(t)]$ is an equivalence class of integer polynomials modulo $g(t)$, and $f(C_g)$ is an integer matrix whose Smith form is sought. In the next sections, we derive results describing some features of the Smith canonical forms of $f(C_g)$ in terms of $[f(t)]$ and $g(t)$.

\section{On the Smith form of $f(C_g)$}\label{sec:SmithFormf(C_g)}

If the base ring $R$ is an integral domain, it can be embedded in a closed field~$\mathbb{F}$, namely, the algebraic closure of the field of fractions of $R$. Hence, the matrix $C_g$ has $n$ eigenvalues (counted with multiplicities) in $\mathbb{F}$. In particular, these eigenvalues are the roots of $g(t)$. Moreover, it is well known~\cite{BP94,GohLR09, NN16,NNT17,NP15} that the eigenvectors of the companion matrix $C_g$ associated with an eigenvalue $\theta$ have the form, up to a nonzero constant,
$$ v_\theta =\begin{bmatrix}
\theta^{n-1}& \cdots& \theta & 1
\end{bmatrix}^T \in \mathbb{F}^n, \qquad \theta : g(\theta)=0.$$
If we assume that $g(t)$ has $n$ distinct roots, then this implies that  $C_g$ is sent to its Jordan canonical form (over $\mathbb{F}$) via similarity by a Vandermonde matrix. If $g(t)$ has multiple roots, the similarity matrix is a confluent Vandermonde. For more details on these classical facts see, for example, \cite{BP94,GohLR09,NN16,NNT17} and the references therein. The matrix $f(C_g) = \mathcal{M}([f(t)])$ therefore has eigenpairs of the form $(f(\theta),v_\theta)$, and in particular, its rank $r$ is equal to the number of the roots of $g(t)$ which are not also roots of $f(t)$, counted with multiplicity. Furthermore, this simple argument shows
the determinant formula~(\ref{eq:detfCg}) of the Introduction. Although typically not stated in this generality, this result is well known at least for some popular choices of $R$ and $g(t)$, for example, when $R$ is any subring of $\mathbb{C}$ and $g(t)=t^n-1$ (so $R_g$ is the ring of circulant matrices) \cite{BP94,Davis}.

We now focus on the case where $R$ is an EDD (note that this implies that $R$ is an integral domain, as in the discussion above), with the goal of studying the Smith form of $f(C_g)$. Recall that every EDD is a B\'{e}zout domain, and therefore a GCD domain. This implies that $R[t]$ is also a GCD domain; that is, given any pair of polynomials $f(t), g(t) \in R[t]$ their gcd exists in $R[t]$. In the following, we will use this fact without further justification.

\section{Proving Theorem~\ref{thm:C}}\label{subsec:f(t)andg(t)notcoprime}

In this section we prove Theorem~\ref{thm:C}. The first step is the technical Lemma~\ref{lem:fatherlemma}, which shows that if $a(t),b(t) \in R[t]$ are two monic polynomials yielding the factorization $g(t)=a(t) b(t)$ then $C_g$ is similar over $R$ to another matrix $X_{a,b} \in R^{n \times n}$ which somehow explicitly displays the factorization of $g(t)$; furthermore, the similarity can be expressed via a special matrix $U_a$: a unit of $R^{n \times n}$ that is completely determined by $a(t)$.
\begin{lemma}\label{lem:fatherlemma}
Suppose that $g(t)=a(t) b(t)$ is a polynomial of degree $n$ where $a(t),b(t)$ are two monic polynomials in $R[t]$. Assume that the degree of $a(t)$ is $m$, define $r:=n-m$, and write $$a(t) = t^m + \sum_{i=0}^{m-1} a_i t^i, \qquad b(t) = t^r + \sum_{i=0}^{r-1} b_i t^i.$$ Denote by $C_g$, $C_a$, $C_b$ the companion matrices of the polynomials $g(t), a(t), b(t)$ respectively, and  let $U_a$ be the $n \times n$ unimodular upper triangular Toeplitz matrix
\[ U_a = \begin{bmatrix}
1 & a_{m-1} & \dots & a_0 &  & \\
& 1 & a_{m-1} &\dots & a_0 &  \\
& & \ddots & \ddots & & \vdots \\
& & & 1 & a_{m-1} & a_{m-2}\\
& & & & 1 & a_{m-1}\\
& & & & & 1
\end{bmatrix} \in R^{n \times n}.    \] Then
$$ U_a C_g U_a^{-1} = X_{a,b} : = \begin{bmatrix}
C_b & 0\\
e_1 e_{r}^T & L_m C_a^T L_m \end{bmatrix}$$
where $e_i$ denotes a vector of size coherent with the matrix partition whose $i$-th entry is $1$ and all of its other entries are $0$,  and $L_m$ is the $m\times m$ flip matrix (i.e. the matrix with $1$'s on the antidiagonal and zeroes elsewhere).
\end{lemma}
 Before proving Lemma~\ref{lem:fatherlemma}, we observe that $(C_b, e_{r}^T, e_1)$ and, respectively, $(L_m C_a^T L_m, e_m^T, e_1)$ are standard triples~\cite{GohLR09} for, respectively, $b(t)$ and $a(t)$. It is therefore a consequence of (a minor modification of)~\cite[Theorem~3.1]{GohLR09} that for all $t$, for some matrix $S$ invertible over the field of fractions of $R$
$$g(t)^{-1} = e_n^T (C_g - t I)^{-1} e_1 = e_n^T (S X_{a,b} S^{-1} - t I)^{-1} e_1;$$
see also~\cite[Theorem 2]{C2GS2} for a related result, also stated for $R=\mathbb{C}$ but discussing more general pencils. It follows in particular that $C_g$ and $X_{a,b}$ are similar \emph{over the field of fractions of $R$}. Lemma~\ref{lem:fatherlemma} goes a step further by showing that one can take $S=U_a$, i.e., $U_a C_g = X_{a,b} U_a$; this is crucial in our context, as it implies that $C_g$ and $X_{a,b}$ are actually similar \emph{over $R$}, and hence, equivalent.

\begin{proof}[Proof of Lemma~\ref{lem:fatherlemma}]
Introduce the vectors
\begin{alignat*}{1}
\alpha &= \begin{bmatrix}a_{m-1} & \cdots & a_0 & 0 & \cdots & 0\end{bmatrix}^T \in R^{n-1},\\
\beta  &= \begin{bmatrix}b_{r-1} & \cdots & b_0 & 0 & \cdots & 0\end{bmatrix}^T \in R^{n-1},\\
\gamma  &= \begin{bmatrix}g_{n-1} & \cdots & g_1\end{bmatrix}^T \in R^{n-1}.
\end{alignat*}
Partition first
$$  U_a C_g = \begin{bmatrix} 1 & \alpha^T\\
0& \widehat{U}_a \end{bmatrix} \begin{bmatrix} -\gamma^T & -g_0 \\
I & 0  \end{bmatrix} =  \begin{bmatrix} \alpha^T - \gamma^T & -g_0\\
\widehat{U}_a & 0\end{bmatrix}.$$
On the other hand, partition
$$ X_{a,b} U_a =  \begin{bmatrix}- \beta^T & 0\\
I & -L_{n-1} \alpha \end{bmatrix} \begin{bmatrix}\widehat{U}_a& L_{n-1} \alpha\\
0^T & 1 \end{bmatrix}= \begin{bmatrix} -\beta^T \widehat{U}_a & - \beta^T L_{n-1} \alpha\\
\widehat{U}_a & 0 \end{bmatrix}.$$
Expanding $g(t)=a(t) b(t)$ in the monomial basis gives
$$\begin{bmatrix}
1 & \beta^T
\end{bmatrix} U_a = \begin{bmatrix} 1 & \gamma^T \end{bmatrix},$$and hence,
$$\beta^T \widehat{U}_a = \gamma^T - \alpha^T.$$
Moreover, since by construction the last $m-1$ entries of the vector $\beta$ and the last $r-1$ entries of the vector $\alpha$ are zero, if $i \neq r$ then $\beta_i (L_{n-1} \alpha)_i = 0$. Hence, $\beta^T L_{n-1} \alpha = \beta_r (L_{n-1} \alpha)_r = b_0 a_0 = g_0.$
\end{proof}

\begin{remark}
\em Combining Lemma~\ref{lem:fatherlemma} with Lemma \ref{lem:triangular} below and known divisibility relations between invariant factors of a submatrix and a matrix (see e.g. \cite{Q97}) yields a potentially interesting consequence. Namely, if $b(t)$ is any monic polynomial that divides $g(t)$, one can write down divisibility relations between the invariant factors of $f(C_g)$ and those of $f(C_b)$. This is a useful property when $\deg b(t) \ll \deg g(t)$, as in this situation the size of $f(C_b)$ is much smaller than the size of $f(C_g)$ and so it is easier to compute the invariant factors of $f(C_b)$ than those of $f(C_g)$. A full discussion is beyond the scope of the present paper.\em
\end{remark}

We now exhibit explicitly the Smith form of $a(C_g)$, where $a(t)$ is any monic divisor (in $R[t]$) of $g(t)$ having degree $m$.

\begin{lemma}\label{thm:A}
If $a(t)$ is a monic divisor (in $R[t]$) of $g(t)$ having degree $m=n-r$, then there exists a unimodular $U \in R^{n \times n}$ such that
$$ Ua(C_g) U_a^{-1} = \begin{bmatrix} I_{r} & 0 \\
0 & 0\end{bmatrix}$$
is in Smith form, where $U_a \in R^{n \times n}$ is the unimodular matrix defined in Lemma~\ref{lem:fatherlemma}.
\end{lemma}

\begin{proof}
We can partition
$$U_a = \begin{bmatrix}
U_{11} & U_{12}\\
0 & U_{22}
\end{bmatrix}, \qquad a(C_g) = \begin{bmatrix}
A & B\\
U_{11} & U_{12}
\end{bmatrix}$$
where $U_{11} \in R^{r \times r}, U_{12} \in R^{r \times m},  A \in R^{m \times r}, U_{22}, B \in R^{m \times m}$. Construct the unimodular matrix
$$ U = \begin{bmatrix}
0 & I_{r}\\
I_m & 0
\end{bmatrix} \begin{bmatrix}
I_m & -A U_{11}^{-1}\\
0 & I_{r}
\end{bmatrix} = \begin{bmatrix}
0 & I_{r}\\
I_m & -AU_{11}^{-1}
\end{bmatrix}$$
so that
$$U a(C_g) U_a^{-1} = \begin{bmatrix}
U_{11} & U_{12}\\
0 & B-AU_{11}^{-1} U_{12}
\end{bmatrix} \begin{bmatrix}
U_{11}^{-1} & -U_{11}^{-1} U_{12} U_{22}^{-1}\\
0 & U_{22}^{-1}
\end{bmatrix} = \begin{bmatrix} I_{r} & 0 \\
0 & K\end{bmatrix},$$
where $K=(B-A U_{11}^{-1} U_{12}) U_{22}^{-1}.$ However, the invertibility of $U$ and $U_a$ implies that $r + \rank K = \rank a(C_g) = r,$ and hence, $K=0$.
\end{proof}

The next technical lemma is useful to reduce the amount of explicit matrix calculations in other proofs; it is well known in matrix theory at least for the case where $R$ is a field \cite[Theorem 1.13(f)]{Higham}, and it can be proved similarly for a general $R$.
\begin{lemma}\label{lem:triangular}
Suppose that $f(t) \in R[t]$ and $A,B$ are square matrices. Then
$$T=\begin{bmatrix}
A & 0\\
X & B
\end{bmatrix} \in R^{n \times n} \Rightarrow f(T) = \begin{bmatrix}
f(A) & 0\\
\star & f(B)
\end{bmatrix} \in R^{n \times n}.$$
Here $\star$ denotes a, possibly nonzero, block of the same size as $X$.
\end{lemma}

We now have all the ingredients to prove Theorem~\ref{thm:C}.

\begin{proof}[Proof of Theorem~\ref{thm:C}]
In this proof, we specialize the notation of Lemmata \ref{lem:fatherlemma} and \ref{thm:A} to the choice $a(t)=z(t)$, $b(t)=G(t)$. Then $U z(C_g) U_z^{-1} = I_{r} \oplus 0.$ By Lemma~\ref{lem:fatherlemma}, writing  $F(t) = \sum_{i=0}^{r-1} F_i t^i,$ we have
\begin{alignat*}{1}
U_z F(C_g) U_z^{-1} &= \sum_{i=0}^{r-1} F_i U_z C_g^i U_z^{-1}=  \sum_{i=0}^{r-1} F_i X_{z,G}^i\\
&\qquad = F(X_{z,G})=
\begin{bmatrix} F(C_G)&0\\
\star & F(L_m C_z^T L_m)\end{bmatrix}
\end{alignat*}
where we used Lemma \ref{lem:triangular} and $\star$ denotes a block element whose precise nature is unimportant. Hence, using Lemma \ref{thm:A},
$$ U f(C_g) U_z^{-1} = U z(C_g) U_z^{-1} U_z F(C_g) U_z^{-1} = F(C_G) \oplus 0. $$
\end{proof}

\section{Factorizing $f(t),g(t)$}\label{subsec:factorisingfg}

In this section we consider factors of $f(t)$ and of $g(t)$. Our first result considers a factorization $f(t)=f_1(t)f_2(t)$ of $f$ and relates the Smith form of $f(C_g)$ to the Smith forms of $f_1(C_g)$, $f_2(C_g)$.  It is known \cite[Theorem II.15]{Newmanbook} that if $A$ and $B$ have coprime determinants then the Smith form of $AB$ is the product of the Smith forms of $A$ and $B$. This immediately proves the following theorem as a special case.

\begin{theorem}\label{thm:factorfbetter}
Let $f(t)=f_1(t) f_2(t)$ and suppose that $\res(f_1,g)$ and $\res(f_2,g)$ are coprime. Denote by $S,S_1,S_2$ the Smith forms of, respectively, $f(C_g)$, $f_1(C_g)$, $f_2(C_g)$. Then $S=S_1 S_2.$
\end{theorem}

The following result is a corollary of Theorem~\ref{thm:factorfbetter}; however, we provide a more elementary proof that does not rely on this theorem.

\begin{corollary}\label{thm:factorf}
Let $f(t)=f_1(t) f_2(t)$ and $\res(f_2,g)$ is a unit of $R$. Then $f_1(C_g) \sim f(C_g)$.
\end{corollary}

\begin{proof}
By~(\ref{eq:detfCg})  the determinant $\mathrm{det}(f_2(C_g))$ is a unit so $f_2(C_g)$ is unimodular. Hence
\(f(C_g) \sim f(C_g) f_2(C_g)^{-1} = f_1(C_g).\)
\end{proof}

Our next result considers a factorization $g(t)=g_1(t)g_2(t)$ of $g(t)$ and relates the matrix $f(C_g)$ to the matrices $f(C_{g_1})$,$f(C_{g_2})$.

\begin{theorem}\label{thm:factorgbetter} Let $g(t)=g_1(t) g_2(t)$ and suppose that $\res(f,g_1)$ and $\res(f,g_2)$ are coprime. Then $f(C_g) \sim f(C_{g_1}) \oplus f(C_{g_2})$.
\end{theorem}

\begin{proof}
It follows from Lemma~\ref{lem:fatherlemma} and Lemma \ref{lem:triangular} that
\[f(C_g) \sim \begin{bmatrix} f(C_{g_2}) & 0\\
X & L_m f(C_{g_1})^T L_m\end{bmatrix}.\]
By \cite[Lemma 6.11]{M4skew} (which is stated for $R=\mathbb{F}[x]$, but in fact only relies on the existence of the Smith form and on Bezout's identity, both valid over every EDD), if the determinants of $A$ and $B$ are coprime then
\[\begin{bmatrix}
A & 0\\
C & B\end{bmatrix} \sim \begin{bmatrix}
A & 0\\
0 & B\end{bmatrix}.\]
To conclude the proof note that
\(L_m f(C_{g_1})^T L_m \sim f(C_{g_1})^T \sim f(C_{g_1}).\)
\end{proof}

\begin{corollary}\label{thm:factorg}
Let $g(t) = g_1(t) g_2(t)$ and $\res(f,g_2)$ is a unit of $R$. Then $f(C_g) \sim I_{\deg (g_2(t))} \oplus f(C_{g_1}).$
\end{corollary}

\begin{proof}
Since $\res(f,g_2)$ is a unit, (\ref{eq:detfCg}) implies that $f(C_{g_1})$ is unimodular and hence is equivalent to the identity matrix. The result then follows from Theorem \ref{thm:factorgbetter}.
\end{proof}

\section{Application to cyclically presented groups and Brieskorn manifolds}\label{sec:applicationtoDunwoody}

The polynomial $g(t)=t^n-1$ and Alexander polynomial $f(t)$ of the torus knot $K(r,s)$ -- see~(\ref{eq:AlexanderTorus}) -- can each be written as a product of cyclotomic polynomials. Before we calculate the Smith form of $f(C_g)$ in Theorem~\ref{thm:dunwoody}, we first calculate the Smith form of $\Phi_m(C_{\Phi_n})$ in Theorem~\ref{thm:smithofcyclofcycl}: this simpler case is potentially useful per se, and it serves the purpose of illustrating some of the basic ideas that we will also use later.
Theorem~\ref{thm:smithofcyclofcycl} can be viewed as a considerable generalization of~\cite[Theorems~2 and 3]{Apostol} which assert that if $m\geq n\geq 1$ then the resultant $\res(\Phi_m,\Phi_n)$ is zero if $m=n$, is $p^{\phi(n)}$ if $m=np^k$ where $k\geq 1$ and $p$ is prime, and is $1$ otherwise. In~\cite[Theorem~3]{Cremona} one step further is taken to derive an expression for $\res(\Phi_m,t^n-1)$. We make repeated use of both these resultant formulae in the proofs of Theorem~\ref{thm:smithofcyclofcycl} and Theorem~\ref{thm:dunwoody}.

The first step is to characterize the first determinantal divisor which we do in Lemma \ref{lem:gamma1}. Recall that the \emph{content} of a polynomial $f(t)=\sum_{i=0}^m f_i t^i \in R[t]$, where $R$ is a GCD domain, is $\mathrm{cont}(f) = (f_0,f_1,\dots,f_m).$

\begin{lemma}\label{lem:gamma1}
Let $g(t) \in R[t]$ be monic of degree $n$, and let $f(t) \in R[t]$ where $R$ is a GCD domain; moreover let $h(t)$ be the unique polynomial of degree less than $n$ such that $f(t) \equiv h(t)\bmod g(t)$.  Then the first determinantal divisor of $f(C_g)$ is (up to units of $R$) $\gamma_1 = \mathrm{cont}(h)$.
\end{lemma}
\begin{proof}
Write $h(t)=\sum_{k=0}^{n-1} h_k t^k$. Recalling Theorem \ref{thm:zero} and the remarks before it, we have $f(C_g)=h(C_g)$. It is readily verified, by finite induction, that the bottom row of $C_g^k$ is equal to $e_{n-k}^T$ for all $k=0,\dots,n-1$ (see also the remarks after Theorem \ref{thm:zero}). It follows that the bottom row of $h(C_g)$ contains as entries precisely the $h_k$, and hence, $\gamma_1 \mid \mathrm{cont}(h).$ On the other hand,
$$(h(C_g))_{ij} = \sum_{k=0}^{n-1} h_k (C_g^k)_{ij},$$
and therefore \emph{all} the entries of $h(C_g)$ are $R$-linear combinations of the coefficients of $h(t)$. This implies $\mathrm{cont}(h) \mid \gamma_1$, and concludes the proof.
\end{proof}

The next steps are Lemmata \ref{lem:pdividesthedifference} and \ref{lem:congruence}: two simple properties of polynomials that will also be handy in proving Theorem \ref{thm:dunwoody}. Lemma  \ref{lem:pdividesthedifference} is a simple consequence of the fact that the $p$-th power map is a ring homomorphism on integers modulo $p$.
\begin{lemma}\label{lem:pdividesthedifference}
Let $p \in \Z$ be a prime and $f(t) \in \Z[t]$ a polynomial. Then for all $k \geq 1$, $p$ divides $[f(t)]^{p^k} - f(t^{p^k})$.
\end{lemma}

\begin{proof}
The proof is by induction on $k$.

When $k=1$, write $f(t) = \sum_{i=0}^{F} f_i t^{i}$, where $F=\deg f(t)$. Applying the multinomial theorem and then Fermat's little theorem in turn, working mod~$p$, we have
\[ f(t)^p\equiv \sum_{i=0}^F (f_it^i)^p =\sum_{i=0}^F f_i^p t^{ip} \equiv \sum_{i=0}^F f_i t^{ip} =f(t^p). \]

Now assume that, for $\ell \in \{ 1,k-1 \}$ and for any $f(t) \in \Z[t]$, $[f(t)]^{p^l}\equiv f(t^{p^l})$~mod~$p$. Set $q=p^{k-1}$ so that $p^k=pq$. Applying the inductive assumption twice, first (for $\ell=k-1$) to $g(t)=[f(t)]^p$ and then (for $\ell=1$) to $h(t)=f(t^q)$, we see that
$ [f(t)]^{pq} \equiv [f(t^q)]^p \equiv f(t^{pq})\bmod p.$
\end{proof}

\begin{lemma}\label{lem:congruence}
Let $R$ be a GCD domain, $g(t) \in R[t]$ be monic, and $f(t) \equiv h(t)\bmod g(t)$ with $\deg h(t) < \deg g(t)$. Then $\mathrm{cont}(f) \mid \mathrm{cont} (h)$.
\end{lemma}

\begin{proof}
Let $m,r,n$ be the degrees of $f(t)=\sum_{i=0}^m f_i t^i$, $h(t)=\sum_{i=0}^r h_i t^i$, $g(t)=\sum_{i=0}^{n} g_i t^i$ (with $g_n=1$) respectively.
Without loss of generality we can suppose that (1) $\mathrm{cont}(f)$ is not a unit of $R$, (2) $r>0$ and (3) $r < m$, else, respectively, (1) $\mathrm{cont}(h) = \mathrm{cont}(f) [ \mathrm{cont}(f)^{-1} \mathrm{cont}(h)]$, (2) $h(t)=0$ and (3) $h(t)=f(t)$, and in all cases the statement becomes obvious.
By assumption there exists $q(t) \in R[t]$ such that $f(t)=h(t)+g(t)q(t)$. Equivalently, there is $\Theta \in R^{m-r+1}$ such that $\Phi = \Gamma + E \Theta$ where
\begin{alignat*}{1}
\Phi = \begin{bmatrix}
f_m & \cdots &f_0
\end{bmatrix}^T&, \qquad
\Gamma = \begin{bmatrix}
0& \cdots & 0 & h_r & \cdots & h_0
\end{bmatrix}^T \in R^{m+1},\\
E&= \begin{bmatrix}
g_n &  & & \\
\vdots & \ddots & & \\
g_1 & & \ddots & \\
g_0 & \ddots &  & g_n \\
& \ddots & \ddots & \vdots\\
& & \ddots & g_1\\
& & & g_0
\end{bmatrix} \in R^{(m+1) \times (m-r+1)}.
\end{alignat*}
Let $T$ be the square submatrix of $E$ obtained by keeping the top $m-r+1$ rows. Manifestly, $\det T=(g_n)^{m-r+1}=1$ so $T$ is unimodular. Hence, noting that the top $m-r+1$ entries of $\Gamma$ are all zero and that $\mathrm{cont}(f)$ divides $(f_n,f_{n+1},\dots,f_m)$,
$$ \Theta = T^{-1} \begin{bmatrix}
f_m & \cdots & f_{n}
\end{bmatrix}^T \Rightarrow q(t) \equiv 0\bmod \mathrm{cont}(f).$$
It follows that $h(t) = f(t) - g(t) q(t)  \equiv 0\bmod \mathrm{cont}(f)$, as required.
\end{proof}

\begin{theorem}\label{thm:smithofcyclofcycl}
Let $m \geq n\geq 1$ and let $\Phi_m(t),\Phi_n(t)\in \Z[t]$ be cyclotomic polynomials. Then the Smith form of $\Phi_m(C_{\Phi_n})$ is
$$S=\begin{cases}
0 \quad &\mathrm{if} \ m=n;\\
p I\quad &\mathrm{if} \ m=n p^k,~\mathrm{where}~p~\mathrm{is~prime}~\mathrm{and}~k\geq 1;\\
I\quad &\mathrm{otherwise}.
\end{cases}$$
\end{theorem}

\begin{proof}
Throughout this proof we use Theorem~\ref{thm:zero} extensively without further reference. The first line follows by the Cayley-Hamilton Theorem; in the third line by (\ref{eq:detfCg}) the determinant $|\mathrm{det} \Phi_m(C_{\Phi(n)})|=|\res(\Phi_m,\Phi_n)|=1$ so $S=I$; and when $n=1$ or $2$ we have $C_{\Phi(n)}=(-1)^{n+1}$, giving the result. Thus we assume that $m=n p^k$, $n\geq 3$,  $p$ prime. We split the proof according to three cases: either (1) $n=p^h$ for some $h > 0$ or (2) $(n,p)=1$ or (3) $n=p^h \ell$ where $h>0$, $\ell > 1$, and $(\ell,p)=1$.

\noindent Case 1 : $n=p^h$. Denote $q=p^{h-1}$ and $r=p^k$. It follows from~\cite[Exercise 12, p.\,237]{Cox} that $\Phi_{p^h}(t)=\Phi_p(t^{p^{h-1}})$ and so  $\Phi_{n} (t) = [t^{pq}-1][t^{q}-1]^{-1}.$ Therefore $1 \equiv t^{pq}$ mod $\Phi_n(t)$; for all $j \in \mathbb{N}$, since $t^{pq}-1$ divides $t^{j qr}-1$, this implies $t^{jqr} \equiv 1$ mod $\Phi_n(t)$. We conclude that $\Phi_m(t) = \sum_{j=0}^{p-1} t^{j qr}\equiv p$~mod~$\Phi_n(t)$ and hence $\Phi_m(C_{\Phi_n})=p I.$

\noindent Case 2 : $(n,p)=1$.
Denote $q=p^{k-1}$. Then $\Phi_m(t)=[\Phi_n(t^{pq})][\Phi_n(t^{q})]^{-1}$ \cite[p.\,160]{Nagell}.
Specializing Lemma \ref{lem:pdividesthedifference} to the polynomial $\Phi_n(t^q)$, we see that $p$ divides $[\Phi_n(t^q)]^p - \Phi_n(t^{pq})$. Dividing by $\Phi_n(t^q)$ and taking into account that $\Phi_n(t)$ divides $\Phi_n(t^q)$, this implies the existence of $\Psi(t) \in \Z[t]$ such that $ \Phi_m(t) \equiv p \Psi(t)$ mod $\Phi_n(t)$. By Lemma \ref{lem:gamma1}, the first determinantal divisor of $\Phi_m(C_{\Phi_n})$, say $\gamma_1$, is the content of the unique polynomial, $\eta(t)$ (say),  \emph{of degree less than $\phi(n)$} and equivalent to $\Phi_m(t)$~ mod~$\Phi_n(t)$. Since $\eta(t)\equiv p\Psi(t)$~mod~$\Phi_n(t)$ Lemma \ref{lem:congruence} implies that $p \mid \mathrm{cont}(p\Psi) \mid \gamma_1$. This fact, together with $\gamma_n=$ $|\mathrm{det} \Phi_m(C_{\Phi_n})|=|\mathrm{Res}(\Phi_m,\Phi_n)|=$ $p^{\phi(n)}$, yields the statement.

\noindent Case 3 : $n=p^h \ell$ where $h>0$, $\ell > 1$, and $(\ell,p)=1$.
Denote $r=p^k$. Then $\Phi_m(t)=\Phi_n(t^r)$ \cite[p.\,160]{Nagell}. By Lemma \ref{lem:pdividesthedifference}, $p$ divides $[\Phi_n(t)]^r-\Phi_m(t)$. The same argument as Case 2 can then be used.
\end{proof}

Now we can prove Theorem~\ref{thm:dunwoody}.

\begin{proof}[Proof of Theorem~\ref{thm:dunwoody}]
It is convenient to split the proof in two cases: (1) $x=1$ (2) $x>1$.

\noindent Case 1 : $x=1$. We have
$$g(t)=\prod_{\delta \in \mathcal{G}} \Phi_\delta(t), \qquad f(t)=\prod_{d \in \mathcal{F}} \Phi_d(t)$$
with
$$\mathcal{G}=\{ \delta \mid n\}, \qquad \mathcal{F}=\{d\mid rs, d\nmid r, d \nmid s\}.$$
Note that $\mathcal{G}\cap\mathcal{F}=\emptyset$ so $(f(t),g(t))=1$ and hence the Smith form for $f(C_g)$ has no zero invariant factors. Let $r=p_1^{\alpha_1} \cdots p_\ell^{\alpha_\ell}$ be the prime factorization of $r$ and $d \in \mathcal{F}$. Then $|\res(g,\Phi_d)|=1$ unless $D=d/(d,n)$ is a positive prime power. This can only happen if $D$ is a positive power of a prime factor $p_i$ of $r$,  for if $D$ divides $s$ then, since $(r,n)=1$, $d=D(d,n)$ is coprime with $r$ so $d|rs$ implies $d|s$, a contradiction. In turn, this implies that $d=p_i^{\beta} k$ with $1 \leq \beta \leq \alpha_i$ and $1 \neq k\mid y.$ Indeed, it cannot be $k=1$, otherwise $d$ divides $r$. Hence, in view of Theorem \ref{thm:factorfbetter}, the sought Smith form is the product of the Smith forms of $f_i(C_g)$, $i=1,\dots,\ell$, with
$$f_i(t) = \prod_{d \in \mathcal{F}_i} \Phi_d(t), \qquad \mathcal{F}_i =\{ p_i^{\beta} k\ :\ 1\leq \beta\leq \alpha, \neq k|y\} \subseteq \mathcal{F}.$$
Moreover, by Corollary~\ref{thm:factorg}, $f_i(C_g) \sim I \oplus f_i(C_h)$ where $h(t)$ is the product of cyclotomics over
$$\mathcal{G}'=\{ 1 \neq \delta\mid y \}.$$
Hence, $f_i(C_h)$ has size precisely $\sum_{1 \neq k\mid y}\phi(k)=y-1$. Furthermore, its determinant is in absolute value
$$|\det f_i(C_h)| = \prod_{d \in \mathcal{F}_i} |\res(\Phi_d,h)| =\prod_{\beta=1}^{\alpha_i} \prod_{1 \neq k \mid y} p_i^{\phi(k)} = \prod_{\beta=1}^{\alpha_i} p_i^{y-1} = \left(p_i^{\alpha_i}\right)^{y-1}.$$ We now claim that $f_i(t) \equiv \Psi_i(t)\bmod h(t)$, with $p_i^{\alpha_i} \mid \Psi_i(t)$. This implies, following an argument analogous to that of Case 2 in the proof of Theorem \ref{thm:smithofcyclofcycl}, that $p_i^{\alpha_i}$ divides the first invariant factor of $f_i(C_h)$, and hence, $f_i(C_h) \sim p_i^{\alpha_i} I_{y-1}$. Since we can repeat this argument for all $i=1, \dots, \ell$, we conclude by Theorem \ref{thm:factorfbetter} that
$$f(C_g) \sim I_{n+1-y} \oplus r I_{y-1}.$$
We now prove the claim. Observe that if $d=p_i^{\beta} k$ then letting $q_{i \beta}=p_i^{\beta-1}$ we have
$$\Phi_d(t) = \frac{\Phi_{k}(t^{p_iq_{i\beta}})}{\Phi_{k}(t^{q_{i\beta}})} \Rightarrow f_i(t) = \prod_{\beta=1}^{\alpha_i}  f_{i\beta}(t), \quad f_{i\beta}(t):=    \frac{h(t^{p_i^{\beta}})}{h(t^{q_{i\beta}})}.$$ By Lemma \ref{lem:pdividesthedifference}, $p_i$ divides $[h(u)]^{p_i} - h(u^{p_i})$ for any variable $u$: dividing by $h(u)$, this in turn implies that $p_i$ also divides the polynomial $[h(u)]^{p_i-1} - h(u^{p_i})/h(u)$. This in particular holds when $u=t^{q_{i \beta}}$, for all $ 1\leq \beta \leq \alpha_i$. In this case $q_{i \beta}$ is a prime power and coprime with $y$. Hence, for all $\delta \mid y$, $\Phi_\delta(t)$ divides $
\Phi_\delta(u)$, and thus $h(t)$ divides $h(u)$. Therefore, for all $\beta=1,\dots,\alpha_i$, there exists a polynomial $\Psi_{i \beta}(t)$ such that $f_{i\beta}(t) \equiv p_i \Psi_{i \beta}(t)\bmod h(t)$. Now let $\Psi_i(t) := p_i^{\alpha_i} \prod_{\beta=1}^{\alpha_i} \Psi_{i \beta}(t)$. Manifestly $p_i^{\alpha_i}$ divides $\Psi_i(t)$, and by the above remarks it follows that
\[ f_i(t) \equiv \prod_{\beta=1}^{\alpha_i} p \Psi_{i \beta}(t) \equiv \Psi_i(t)\bmod h(t); \]
this proves the claim.

\noindent Case 2 : $x>1$. Here $z(t)=(f(t),g(t))=\prod_{d\in\Sigma} \Phi_d(t)$ where $\Sigma$ consists of all divisors of $(n,rs)$ that are neither divisors of $(n,r)$ nor $(n,s)$. It follows that $\mathrm{deg}(z)=\sum_{d\in\Sigma} \phi(d)=(x-1)(y-1)$ and so by Lemma \ref{thm:A} the Smith form has $(x-1)(y-1)$ zero invariant factors.  Since $(r,s)=1$ we have $(n,rs)=xy$ so after having removed common factors, as well as the trivial factor $t-1$ in $g(t)$, we are left with the index sets
$$\mathcal{G}=\{1 \neq \delta \mid x\} \cup \{1 \neq \delta \mid y\}\cup \{ \delta \mid n, \delta \nmid xy \},$$
$$\mathcal{F}=\{d\mid sr, d\nmid r, d\nmid s, d\nmid xy \}.$$
Let $d \in \mathcal{F}$ and suppose $\delta \mid n$ but $\delta \nmid xy$. Then $\delta/d$ cannot be a positive or negative prime power, as otherwise $d \mid xy$ or $\delta \mid xy$, respectively; thus $|\res(\Phi_d,\Phi_\delta)|=1$ and by Theorem \ref{thm:factorg} we can effectively (up to neglecting some trivial invariant factors) replace $g(t)$ with $h(t)$, the product of cyclotomics over the set
$$\{1 \neq \delta \mid x\} \cup \{1 \neq \delta \mid y\} =: \mathcal{G}_1 \cup \mathcal{G}_2.$$
Suppose $\delta \in \mathcal{G}_1$ and $d \in \mathcal{F}$. If $\delta\mid x$ then, since $d \nmid r$ and $d \nmid xy$, the only possibility for $\Phi_\delta(t)$ and $\Phi_d(t)$ to have a nontrivial resultant is for $d$ to be of the form $\delta \hat{s}$ where $\hat{s}$ is the power of a prime factor of $s$ and $\hat{s}|s$, $\hat{s}\nmid y$. A similar argument holds if $\delta \in \mathcal{G}_2$, so  we can replace $\mathcal{F}$ with $\mathcal{F}_1 \cup \mathcal{F}_2$ with
$\mathcal{F}_1 := \{\delta \hat{s}: 1 \neq \delta \mid x, \hat{s} $ prime power, $\hat{s} \mid s$, $\hat{s} \nmid y\}$ and $\mathcal{F}_2:= \{\delta \hat{r}:  1 \neq \delta \mid y, \hat{r}$ prime power, $\hat{r} \mid r$, $\hat{r} \nmid x\}.$
Moreover, observing that $s$ and $r$ are coprime, if $d\in\mathcal{F}_i, \delta \in \mathcal{G}_j$ ($\{i,j\}=\{1,2\}$) then $|\res (\Phi_d,\Phi_\delta)|=1$. Thus invoking Theorems~\ref{thm:factorfbetter} and~\ref{thm:factorgbetter}, we also see that the Smith form of $f(C_h)$ is the product of the Smith forms of $I \oplus f_1(C_{h_1})$ and $I \oplus f_2(C_{h_2})$, where the sizes of the identity matrices is clear from the context and, for $i=1,2$, $f_i(t)$ and $h_i(t)$ are products of cyclotomics whose indices vary in $\mathcal{F}_i$ and $\mathcal{G}_i$, respectively.
We now can, in essence, follow the first part of this proof to show that the Smith form of $f_1(C_{h_1})$ has non-unit invariant factors $s/y$ ($x-1$ times), and the Smith form of $f_2(C_{h_2})$ has non-unit invariant factors $r/x$ ($y-1$ times). More precisely, a slight modification is needed to take into account that, when writing  (for example) $f_i(t)=\prod_\beta f_{i \beta}(t)$, the exponents $\beta$ no longer vary between $1$ and $\alpha_i$ but between $\gamma_i+1$ and $\alpha_i$, where $\gamma_i$ is the power of $p_i$ in the prime factorization of $x$. This is a consequence of the fact that, as in the definition of $\mathcal{F}_1$, we must select $\hat{s}$ as a prime power dividing $s$, but not $y$. However, as apart from this subtlety the argument is completely analogous, we omit the details.

Finally, the statement follows by multiplying two diagonal matrices.
\end{proof}

\noindent  \textbf{Acknowledgements} The authors thank Alessia Cattabriga for insightful discussions on Brieskorn and Dunwoody manifolds, Steve Mackey for enlightening comments on the history of the Smith Theorem, and the referees for providing feedback on the presentation. Vanni Noferini  acknowledges support by an Academy of Finland grant (Suomen Akatemian p\"{a}\"{a}t\"{o}s 331240). Gerald Williams was supported for part of this project by Leverhulme Trust Research Project Grant RPG-2017-334.

\end{document}